

\input amstex
\documentstyle{amsppt}
\magnification=1200
\catcode`\@=11
\redefine\logo@{}
\catcode`\@=13

\define \bn{\Bbb N}
\define \bz{\Bbb Z}
\define \bq{\Bbb Q}
\define \br{\Bbb R}
\define \bc{\Bbb C}

\define \M{{\Cal M}}
\define\Ha{{\Cal H}}
\define\La{{\Cal L}}
\define\Da{{\Cal D}}

\define\0o{{\overline 0}}
\define\1o{{\overline 1}}
\define\rk{\text{rk}~}

\TagsOnRight

\topmatter

\title
The remark on Discriminants of K3 Surfaces Moduli as Sets of Zeros
of Automorphic forms
\endtitle

\author
Viacheslav V. Nikulin \footnote{Partially supported by
Grant of Russian Fund of Fundamental Research.
\hfill\hfill}
\endauthor

\address
Steklov Mathematical Institute,
ul. Vavilova 42, Moscow 117966, GSP-1, Russia.
\endaddress
\email
slava\@nikulin.mian.su
\endemail

\abstract
We show that for any $N>0$ there exists a natural even $n>N$ such that
the discriminant of moduli of K3 surfaces of the degree $n$
is not equal to the set of zeros of any automorphic form on
the corresponding IV type domain.

We give the necessary condition on
a "condition $S\subset L_{K3}$ on Picard lattice of K3'' for
the corresponding moduli $\M_{S\subset L_{K3}}$ of K3 to have
the discriminant which is equal to the set of zeros of an automorphic
form. We conjecture that the set  of $S\subset L_{K3}$ satisfying
this necessary condition is finite if $\rk S \le 17$.
We consider this
finiteness conjecture as "mirror symmetric'' to the known
finiteness results for arithmetic reflection groups in hyperbolic
spaces and as important for the theory of Lorentzian Kac--Moody
algebras and the related theory of automorphic forms.
\endabstract

\rightheadtext
{Discriminants of K3 Moduli}

\leftheadtext{Viacheslav V. Nikulin}

\endtopmatter

\document

\subhead
0. Introduction
\endsubhead

The subject of this paper is directly connected with the
recently developed theory of Lorentzian Kac--Moody algebras,
the corresponding theory of automorphic forms on IV type domains,
related geometry of K3 surfaces and their moduli, and mirror symmetry.
See
\cite{Bo1}, \cite{Bo2}, \cite{Bo3}, \cite{G1},\cite{G2},
\cite{G3}, \cite{G4}, \cite{GN1}, \cite{GN2},
\cite{GN3}, \cite{N9}, \cite{N10}.
We refer a reader
who is interested in physical applications to the recent papers
\cite{H--M} and \cite{Ka}.

Results of this article were inspired by my thinking over
recent announcements by J. Jorgensen and A. Todorov on discriminants
of K3 surfaces.
I am grateful to Prof. V.A. Gritsenko for
very useful discussions on this subject.

\subhead
1. Definitions
\endsubhead

Let $L_{K3}$ be an even hyperbolic unimodular lattice of
the signature $(3,19)$. It is known that $L_{K3}$ is isomorphic to
the $H^2(X; \bz)$ for a K3 surface $X$ over $\bc$. To be shorter,
we set $L=:L_{K3}$ below.

Let $S \subset L$ be a primitive hyperbolic sublattice, i.e.
$S$ is a lattice of the signature $(1,k)$ and $L/S$ is a free
$\bz$-module. The pair $S\subset L$ is called {\it" a condition
on Picard lattice of K3''} (or, shortly {\it "condition
$S \subset L$''}). We correspond to the condition $S\subset L$
a {\it lattice $T(S\subset L)=S^\perp_L$ of the signature
$(2,20-\rk S)$} and a {\it complex symmetric domain of the type IV
$$
\Omega(S \subset L)=
\{ \bc \omega \subset T(S\subset L) \otimes \bc\ |
\ \omega^2=0,\ \omega \cdot \overline{\omega}>0 \}
\tag1.1
$$
of the dimension} $20-\rk S$.
Any $e \in S\otimes \br$ with $e^2<0$ defines
a codimension one {\it symmetric
subdomain of the type IV}
$$
D_e=\{ \bc \omega \in \Omega(S\subset L)\ | \ \omega \cdot e=0\}.
\tag{1.2}
$$
The subset
$$
\widetilde{\Cal D}(S\subset L)=
\bigcup_{\delta \in T(S\subset L),\  \delta^2=-2}{D_\delta}
\tag{1.3}
$$
is called the {\it discriminant of the condition $S\subset L$}.

We will consider automorphic forms (i.e., holomorphic, with
some non-negative weight) on
$\Omega (S\subset L)$
with respect to subgroups of finite index
$$
G\subset O(T(S\subset L)).
\tag{1.4}
$$
If there exists such an automorphic form $\Phi$ of positive weight
with the set of zeros $\widetilde{\Cal D}(S\subset L)$, we say
that {\it the discriminant
$\widetilde{\Cal D}(S\subset L)$
is equal to the set of zeros of an automorphic form}.
We remark that
for this definition, we don't fix the subgroup
$G$ of finite index.

We want to show that this situation is extremely rare if
$\rk S\le 17$.

\vskip10pt

The situation which is described above is purely arithmetic and
one can forget about $S\subset L$ considering only
the lattice $T=T(S\subset L)$. Considering of
$S\subset L$ is
important for the theory of K3 surfaces.
One can correspond to
the condition on Picard lattice of K3 $S\subset L$ a family
$\pi: \overline{\Cal X}_{S\subset L} \to
\overline{\M}_{S\subset L}$ of K3 surfaces where
$\overline{\M}_{S\subset L}= G\setminus \Omega(S\subset L)$.
For this family, $X_m=\pi^{-1}(m)$ is a non-singular K3 surface if
$m \in \M_{S\subset L}=
G \setminus (\Omega (S\subset L)-
\widetilde{\Cal D}(S\subset L))$, and $\pi^{-1}(m)$ is
a singular K3 surface with Du Val singularities if
$m \in \Da_{S\subset L}=
G\setminus \widetilde{\Da }(S\subset L)$.
For a general point
$m \in \M_{S\subset L}$, the embedding $S_{X_m}\subset H^2(X_m;\bz)$
of Picard lattice $S_{X_m}$ of $X_m$
is isomorphic to $S\subset L$.
This K3 interpretation follows (see \cite{N1}, \cite{N2})
from general results
\cite{Tjurina, \u S}, \cite{P-\u S--\u S},
\cite{B--R}, \cite{Ku}.
Automorphic forms on $\Omega (S\subset L)$ give sections of
appropriate "linear'' sheafs on
$\overline{\M}_{S\subset L}$. Thus,
geometrically, the discriminant $\widetilde{\Da}(S\subset L)$
is the set of zeros of an automorphic form
if $\Da_{S\subset L}$ is the set of zeros of a section of an
appropriate "linear'' sheaf on $\overline{\M}_{S\subset L}$.
We don't specify the group $G$ since there are several natural
choices of this group. If the discriminant
$\widetilde{\Da }(S\subset L)$ is
equal to the set of zeros of an automorphic form, then for any
subgroup $G\subset O(T(S\subset L))$ of finite index
the preimage of the discriminant $\Da_{S\subset L}$
is the set of zeros of a section of an
appropriate "linear'' sheaf on some finite ramified covering of
$\overline{\M}_{S\subset L}$.

\subhead
2. Results
\endsubhead

Assume that the lattice $S \cong  \langle n \rangle$, i.e.
$S=\bz h$ is one-dimensional and $h^2=n$
where $n$ is an even natural number. Up to isomorphism
(automorphisms of the lattice $L$), the
condition $S\subset L$ is defined by $n$. (It follows, for example,
from the analog of Witt Theorem in \cite{N1}, \cite{N2}.)
We want to prove

\proclaim{Theorem 2.1} For any $N >0$ there exists $n>N$ such that
for $S=\langle n \rangle$ the discriminant
$\widetilde{\Da}(S \subset L)$
is not equal to the set of zeros of any automorphic form.
\endproclaim

\demo{Proof} Let us take a primitive sublattice $T_1\subset L$ such
that $T_1$ has signature $(2, k)$ where $k \ge 2$ and $T_1$ does not
have elements with square $(-2)$. If $k=2$, additionally suppose
that the lattice $T_1$ does not have an isotropic sublattice of the
rank $2$. (It is known \cite{N2} that any
even lattice $K$ of signature $(p,q)$ where $p\le 3$, $q \le 19$ and
$p+q \le 11$ has a primitive embedding to $L$. It follows that
the lattice $T_1$ does exist.) Let $S_1=(T_1)^\perp_L$.
Then $T_1=T(S_1\subset L)$ where $S_1\subset L$ is some condition
on Picard lattice of K3. We denote $\pi_{S_1}$ and $\pi_{T_1}$
orthogonal projections of $L$ on $S_1$ and $T_1$ respectively
(they are defined over $\bq$).

Let us consider {\it the cone}
$$
V(S_1)=\{ x \in S_1 \otimes \br\ |\ x^2>0\},
\tag{2.1}
$$
and its {\it half-cone} $V^+(S_1)$ and the corresponding
{\it hyperbolic space}
$\La (S_1)=V^+(S_1)/\br_{++}$. We consider
$$
\split
&\Delta^{(2)}(S_1\subset L)=\{ \delta \in L\ | \ \delta^2=-2,\\
&\text{either\ } (\pi_{S_1}(\delta))^2< 0\ \text{or\ }
\pi_{S_1}(\delta )=0, \  \text{either\ }
(\pi_{T_1}(\delta))^2< 0\ \text{or\ } \pi_{T_1}(\delta)=0 \}.
\endsplit
\tag{2.2}
$$
Evidently, there exist $a(S_1)\in \bn$ and $b(T_1)\in \bn$
such that
$$
\pi_{S_1}(\delta) \in (1/a(S_1))S_1
\tag{2.3}
$$ and
$$
\pi_{T_1}(\delta) \in (1/b(T_1))T_1
\tag{2.4}
$$
for any  $\delta \in \Delta^{(2)}(S_1\subset L)$.
It follows that
there are some $A(S_1) \in \bn$ and $B(T_1) \in \bn$ such
that
$$
-A(S_1) \le (a(S_1)\pi_{S_1}(\delta))^2\le 0,\ \
(a(S_1)\pi_{S_1}(\delta))^2\in \bz,
\tag{2.5}
$$
and
$$
-B(T_1) \le (b(T_1)\pi_{T_1}(\delta))^2\le 0,\ \
(b(T_1)\pi_{T_1}(\delta))^2 \in \bz
\tag{2.6}
$$
for any $\delta \in \Delta^{(2)}(S_1\subset L)$.
For $e \in S_1\otimes \br$ with $e^2<0$ we denote
$$
\Ha_e=\{\br_{++}x \in \La(S_1)\ |\ x\cdot e=0\}
\tag{2.7}
$$
a {\it hyperplane} in the hyperbolic space
$\La(S_1)$ which is {\it orthogonal to $e$}. Since \thetag{2.3} and
\thetag{2.5}, the set of hyperplanes
$\Ha_{\delta_1}$,
$\delta_1 \in \pi_{S_1}(\Delta^{(2)}(S_1\subset L))-\{0\}$
is locally finite
in the hyperbolic space $\La(S_1)$.

Let  $h \in S_1\cap V^+(S_1)$ be a primitive element of $S_1$.
Since $h^2>0$, the element $h$ defines a sublattice
$[h]\subset S_1\subset L$ of the rank one which gives a condition on
Picard lattice of K3.

We have the following basic

\proclaim{Lemma 2.2} Assume that the lattice $T_1=T(S_1\subset L)$
has $\rk T_1 \ge 4$ and the lattice $T_1$
does not have elements with square $(-2)$. If $\rk T_1=4$, additionally
suppose that the lattice $T_1$ does not have an
isotropic sublattice of the rank $2$. Assume that
for a primitive $h \in S_1 \cap V^+(S_1)$, the condition
$[h]\subset L$ on Picard lattice has the discriminant
$\widetilde{\Da}([h]\subset L)$ which is equal to the
set of zeros of an automorphic form.
Then the point $\br_{++}h$ belongs to a hyperplane
$\Ha_{\delta_1}$,
$\delta_1 \in \pi_{S_1}(\Delta^{(2)}(S_1\subset L))$, in $\La (S_1)$.
\endproclaim

\demo{Proof of Lemma 2.2}
Let us
suppose that the discriminant $\widetilde\Da([h]\subset L)$ is
equal to the set of zeros of an automorphic form $\Phi$ on
$\Omega([h]\subset L)$. It contains symmetric subdomain
$\Omega(S_1\subset L)\subset \Omega([h]\subset L)$.
The restriction $\Phi_1=\Phi | \Omega(S_1\subset L)$ is
an automorphic form on $\Omega (S_1\subset L)$
of the same weight as the weight of $\Phi$. In particular,
it has a positive weight and is not equal to a non-zero constant.
Since $\rk T_1 \ge 4$, the
dimension $\dim \Omega (S_1\subset L)\ge 2$.
By Koecher principle, the automorphic form
$\Phi_1=\Phi | \Omega(S_1 \subset L)$
has a non-empty
set of zeros in $\Omega (S_1\subset L)$. Since the set of zeros of
$\Phi$ in $\Omega ([h]\subset L)$ is equal to the discriminant
$\tilde{\Da} ([h]\subset L)$, it follows that there
exists $\delta \in L$ with the following properties:
$\delta^2=-2$, $h\cdot \delta=0$, either
$\pi_{T_1}(\delta) = 0$ (then $\Phi_1=0$), or
$(\pi_{T_1}(\delta))^2<0$ (then $\Phi_1$ has a zero).
Since
$S_1$ is hyperbolic and $h^2>0$, we then get that
either $(\pi_{S_1}(\delta))^2<0$ or $\pi_{S_1}(\delta)=0$.
The last case is impossible, since $T_1$ does not have
elements $\delta \in T_1$ with $\delta^2=-2$.
It follows that
$\delta\in \Delta^{(2)}(S_1\subset L)$ and
$\br_{++}h \in \Ha_{\pi_{S_1}(\delta)}$ where
$\Ha_{\pi_{S_1}(\delta)}$ is a hyperplane which is orthogonal
to $\pi_{S_1}(\delta)\not= 0$.

It follows the lemma.
\enddemo

Now Theorem 2.1 follows from Lemma 2.2 and the
following simple well-known fact: The set of points $\br_{++}h$,
where $h \in S_1$ and $0<h^2\le N$, is locally finite in $\La(S_1)$.

Since the set of points $\br_{++}h\in \La (S_1)$, $h \in S_1$,
is everywhere dense in $\La (S_1)$ and the set of hyperplanes
$\Ha_{\pi_{S_1}(\delta)}$, $\delta \in \Delta^{(2)}(S_1\subset L)$,
is locally finite in $\La (S_1)$, for any $N>0$ one can find
a point $\br_{++}h \in \La(S_1)$ such that $h\in S_1$ is primitive,
$h^2>N$ and $\br_{++}h$ does not belong to any hyperplane
$\Ha_{\pi_{S_1}(\delta)}$, $\delta \in \Delta^{(2)}(S_1\subset L)$.
By Lemma 2.2, the condition $[h]\subset L$ has the discriminant
which is not equal to the set of zeros of any automorphic form.
\enddemo

Replacing the sublattice $[h]\subset L$ above
by a condition $S\subset L$, similarly, one can prove the following
result:

\proclaim{Theorem 2.3} Let $S\subset L$ be a condition on Picard
lattice of K3 and the discriminant $\widetilde{\Da}(S\subset L)$
is equal to the set of zeros of an automorphic form.

Then for any primitive hyperbolic sublattice
$S_1\subset L$ such that $S\subset S_1$, $\rk S_1 \le 18$,
$T_1=(S_1)^\perp_L$ does not have elements with square $(-2)$ and
if $\rk S_1=18$, the lattice
$T_1$ does not have isotropic sublattices of the rank $2$,
there exists $\delta \in \Delta^{(2)}(S_1\subset L)$
such that $\La (S)\subset \Ha_{\pi_{S_1}(\delta)}$
(equivalently, $S\cdot \delta=0$).
\endproclaim

\demo{Proof} The proof is the same as for Lemma 2.2.
\enddemo

We think that the necessary condition of Theorem 2.3 for
$S\subset L$ to have the discriminant $\widetilde\Da (S\subset L)$
which is equal to the set of zeros of an automorphic form is
very restricted.

\proclaim{Conjecture 2.4} If $\rk S \le 17$,
the set of conditions $S\subset L$
which satisfy the necessary condition of Theorem 2.3 is finite
(up to isomorphism). In particular, if $\rk S \le 17$,
the set of conditions $S\subset L$ with the discriminant
$\widetilde\Da (S\subset L)$ which is equal to the set of zeros
of an automorphic form is finite (up to isomorphism).

If $\rk S=18$ and the lattice $T=T(S\subset L)=S^\perp_L$
contains an isotropic sublattice of the rank one,
then the condition $S\subset L$ has the
discriminant $\widetilde\Da (S\subset L)$ which is equal
to the set of zeros of an automorphic form if and only if
the lattice $T$ either has an element with
square $(-2)$ or an isotropic sublattice of the rank $2$.

If $\rk S=19$ and the lattice $T=T(S\subset L)=S^\perp_L$
contains an isotropic sublattice of the rank one,
the condition $S\subset L$ has the
discriminant $\widetilde\Da (S\subset L)$ which is equal to
the set of zeros of an automorphic form (for this case
$\dim \Omega (S\subset L)=1$).
\endproclaim

Let $S$ be an even hyperbolic lattice and $W^{(2)}(S)\subset O(S)$
the group generated by all reflections in $\delta \in S$ with
$\delta^2=-2$. The lattice $S$ is called {\it $2$-reflective} if
index $[O(S):W^{(2)}(S)]$ is finite.

We think that Conjecture 2.4 is "mirror symmetric" to the
following result:

\proclaim{Theorem 2.5} For $\rk S \ge 3$, the set of $2$-reflective
lattices is finite (\cite{N3}, \cite{N4},
\cite{N5}, \cite{N7}, \cite{N8}).
If $\rk S=2$, the lattice $S$ is $2$-reflective
if and only if $S$ has either an element with square $-2$
or with square $0$
(see \cite{P-\u S--\u S}). If $\rk S=1$, the lattice $S$ is
obviously $2$-reflective.
\endproclaim

We hope to consider Conjecture 2.4 in more advanced publications.

Similarly, one can formulate and prove more general results
replacing elements with square $(-2)$ by elements with squares which
are contained in a fixed finite set of negative integers.
These results are connected
with the general Theory of reflection groups
in hyperbolic spaces, the general Theory of Lorentzian Kac--Moody
algebras and the related Theory of automorphic forms.
See \cite{Bo1}, \cite{Bo2}, \cite{Bo3}, \cite{G1}, \cite{G2},
\cite{G3}, \cite{G4}, \cite{GN1}, \cite{GN2}, \cite{GN3},
\cite{N9}, \cite{N10}.

\enddocument

\end